\documentclass[preprint,12pt]{elsarticle}

\usepackage{amssymb, amsmath, amstext, amsthm, amsfonts, amscd, float, multirow}
\usepackage[utf8]{inputenc}
\usepackage[T1]{fontenc}
\usepackage[mathscr]{eucal}
\usepackage{pseudocode}
\usepackage[dvipsnames]{xcolor}
\usepackage{tikz}

\newcommand{\R}{\mathbb{R}}
\newcommand{\C}{\mathbb{C}}
\newcommand{\W}{\mathbb{W}}
\newcommand{\G}{\mathbb{G}}

\newcommand{\N}{\mathbb{N}}
\newcommand{\abs}[1]{\bigl| #1 \bigr|} 
\newcommand{\norm}[1]{\lVert#1\rVert} 
\newcommand{\Pp}{\mathbb{P}}

\newtheorem{theorem}{Theorem}[section]
\newtheorem{proposition}[theorem]{Proposition}
\newtheorem{remark}{Remark}[section]
\newtheorem{example}{Example}[section]
\newtheorem{definition}{Definition}[section]
\newtheorem{lemma}[theorem]{Lemma}
\newtheorem{corollary}{Corollary}[section]

\theoremstyle{definition}
\newtheorem{experiment}{Experiment}[section]

\journal{Applied Mathematics and Computation}

\begin{document}

\begin{frontmatter}

\title{Isospectral Reductions of Non-negative Matrices}

\author[inst1]{Alexandre Baraviera}

\affiliation[inst1]{organization={Instituto de Matem\'{a}tica e Estat\'{i}stica, UFRGS},
            addressline={Av. Bento Gon\c{c}alves 9500}, 
            city={Porto Alegre},
            postcode={91500}, 
            state={RS},
            country={Brazil}}

\author[inst2]{Pedro Duarte}
\author[inst3]{Longmei Shu}
\author[inst4]{Maria Joana Torres}

\affiliation[inst2]{organization={CMAF, Departamento de Matem\'{a}tica},
            addressline={Faculdade de Ci\^{e}ncias da Universidade de Lisboa}, 
            city={Campo Grande},
            postcode={1749-016}, 
            state={Lisboa},
            country={Portugal}}
\affiliation[inst3]{organization={Department of Mathematical Sciences},
            addressline={Rutgers University - Camden}, 
            city={Camden},
            postcode={08102}, 
            state={NJ},
            country={USA}}
\affiliation[inst4]{organization={CMAT and Departamento de Matem\'atica},
            addressline={Universidade do Minho},
            city={Campus de Gualtar},
            postcode={4700-057},
            state={Braga},
            country={Portugal}
            }

\begin{abstract}
Isospectral reduction is an important tool for network/matrix analysis as it reduces the dimension of a matrix/network while preserving its eigenvalues and eigenvectors. The main contribution of this manuscript is a proposed algorithmic scheme to approximate the stationary measure of a stochastic matrix  based on isospectral reductions. We run numerical experiments that indicate this scheme is advantageous when there is more than one eigenvalue near $1$, precisely the case where iterative methods perform poorly. We give a partial explanation why this scheme should work well, showing that in some situations isospectral reduction improves the spectral gap.
\end{abstract}


\end{frontmatter}

\section{Introduction}

Markov chains are a powerful tool for modeling and predicting the behavior of complex systems. They are used in a wide variety of fields, including finance, biology, and computer science. The stationary distribution of a Markov chain describes its long-term behavior and can be used to understand and control the behavior of Markov chains. One application of the stationary distribution in real world networks like citation networks and the Internet is in computing various network centrality measures, such as PageRank or eigenvector centrality. These measures are often based on the stationary distribution of random walks and help in ranking the importance of nodes in the network (see \cite{austin06, brin98, dorogovtsev03, newman06, noh04}).  In \cite{nowak06}, the author discusses the significance of stationary distributions in understanding the long-term behavior of populations and the stable strategies that emerge through evolutionary processes. Also, in \cite{muller15}, the authors explore the application of stationary distributions to analyze the long-term behavior of populations in various ecological models. Examining the stationary distributions of susceptible and infected individuals in epidemic models highlights their importance for predicting disease dynamics and informing public health interventions.

Typically we compute the stationary vectors for a Markov chain with iterative methods. While these iterative methods work well for sparse transition matrices of a Markov chain with a simple eigenvalue 1, the computation can take a long time as the transition matrices become very large or if there are multiple eigenvalues near 1. 

The classical Perron-Frobenius theorem states that for a non-negative, primitive matrix $A\in \R^{n\times n}$  and a vector $v\in \R^n_+$, the iterates $A^n \, v$ converge geometrically in the projective space to the dominant eigenvector of $A$. Moreover, the logarithmic rate of convergence is $\log r<0$, where $0<r<1$ is the largest absolute value of the ratios of the non-dominant eigenvalues to the dominant eigenvalue.
For more details, see formula (8.2.10) in Horn and Johnson~\cite[Chapter 8]{horn12}. In the case of stochastic primitive matrices, the positive dominant eigenvector corresponds to the stationary distribution, the dominant eigenvalue is 1, and the rate $r$ is given by the largest absolute value of the eigenvalues of $A$ within the open unit disk.
In particular, convergence to the stationary distribution, as guaranteed by the Perron-Frobenius theorem, can be slow when the second-largest eigenvalue of $A$ is close to $1$ in absolute value, i.e., when $r\sim 1$.
Modern Linear Algebra software employs a class of iterative methods to approximate the eigenvectors of large matrices, known as Arnoldi methods, which are variations of a technique proposed by W. E. Arnoldi in~\cite{arnoldi51}. These Arnoldi methods face similar limitations to the iterative scheme suggested by the Perron-Frobenius theorem. In particular, they can be slow and less accurate when one or more eigenvalues are close to the dominant eigenvalue.
See~\cite[Section 6.7]{saad11} for more details.

The recently developed theory of Isospectral Transformations (IT) of matrices and networks allows for advances in various areas and has led to several surprising results \cite{webb14, baraviera22, shu18, bunimovich19, torres14, duarte20,  kempton20, morfonios21, rontgen20,rontgen23, smith19}. The theory of isospectral transformations was initially aimed at reduction (i.e. simplification) of networks,  while keeping all the information about the spectrum of their weighted adjacency, Laplace, or other matrices. However, it turns out that all the information about the eigenvectors of these matrices also gets preserved under ITs \cite{torres14,shu18}. The eigenvectors of the reduced matrix are bijective projections of the eigenvectors of the original matrix and we have formulas to reconstruct the original eigenvectors from the reduced eigenvectors.

Therefore it is natural to ask if we can efficiently use the reductions to compute the stationary measure of a large stochastic matrix in a Markov chain. The main goal of the present paper is to investigate this question. 

The paper is organized as follows. We define isospectral reductions first (Section 2), then introduce various measurements involving the eigenvalues of a stochastic matrix, including the diameter semi-norm which measures the variation of the columns of a stochastic matrix (Section 3). We show that the reduction of a stochastic matrix is still a stochastic matrix and the diameter semi-norm decreases after reduction. We also show that for bi-stochastic matrices the Gershgorin region, which traps the eigenvalues of a matrix by disks around its diagonal entries, shrinks in disk sizes after isospectral reductions. This agrees with the results in \cite{webb12}. A new computational scheme is proposed for the computation of stationary measures for stochastic matrices using isospectral reduction (Section 4). We run numerical experiments to compare this new scheme with traditional methods. Our experiments indicate that the new scheme is faster and more accurate for large stochastic matrices with one or more eigenvalues close to 1. In Section 5, we construct various examples of stochastic matrices that become positive (all entries are positive) after reduction, which correlates to the decrease in the diameter semi-norm. In Section 6 we discuss our results and further research directions.

\section{Isospectral reductions}
In this section we recall definitions of the isospectral reductions of graphs and networks.

Let $\W$ be the set of rational functions of the form $w(\lambda)=p(\lambda)/q(\lambda)$, where $p(\lambda),q(\lambda)\in\C[\lambda]$ are polynomials having no common linear factors, i.e., no common roots, and where $q(\lambda)$ is not identically zero. $\W$ is a field under addition and multiplication \cite{webb14}.

\subsection{Isospectral Graph Reductions}

Let $\G$ be the class of all weighted directed graphs with edge weights in $\W$. More precisely, a graph $G\in\G$ is an ordered triple $G=(V,E,w)$ where $V=\{1,2,\dots,n\}$ is the \emph{vertex set}, $E\subset V\times V$ is the set of \emph{directed edges}, and $w:E\to\W$ is the \emph{weight function}. Denote by $M_G=(w(i,j))_{i,j\in V}$ the \emph{weighted adjacency matrix} of $G$, with the convention that $w(i,j)=0$ whenever $(i,j)\not\in E$. We will alternatively refer to graphs as networks because weighted adjacency matrices can be used to define static (i.e. non evolving) real world networks.

Observe that the entries of $M_G$ are rational functions. Let's write $M_G(\lambda)$ instead of $M_G$ here to emphasize the role of $\lambda$ as a variable. For $M_G(\lambda)\in\W^{n\times n}$, we define the spectrum, or multiset of eigenvalues to be $$\sigma(M_G(\lambda))=\{\lambda\in\C:\det(M_G(\lambda)-\lambda I)=0\}.$$ 
Notice that $\sigma(M_G(\lambda))$ can have more than $n$ elements, some of which can be the same. 

Throughout the rest of the paper, the spectrum is understood to be a set that includes multiplicities.



An eigenvector for an eigenvalue $\lambda_0\in\sigma(M_G(\lambda))$ is defined to be $u\in\C^n,u\neq0$ such that 
$$M_G(\lambda_0)u=\lambda_0u.$$ 
One can see that the eigenvectors of $M_G(\lambda)\in\W^{n\times n}$ for $\lambda_0$ are the same as the eigenvectors of $M_G(\lambda_0)\in\C^{n\times n}$ for $\lambda_0$. Similarly the generalized eigenvectors of $M_G(\lambda)$ for $\lambda_0$ are the generalized eigenvectors of $M_G(\lambda_0)$ for $\lambda_0$.

A path $\gamma=(i_0,\dots,i_p)$ in the graph $G=(V,E,w)$ is an ordered sequence of distinct vertices $i_0,\dots,i_p\in V$ such that $(i_l,i_{l+1})\in E$ for $0\le l\le p-1$. The vertices $i_1,\dots,i_{p-1}\in V$ of $\gamma$ are called \emph{interior vertices}. If $i_0=i_p$ then $\gamma$ is a \emph{cycle}. A cycle is called a \emph{loop} if $p=1$ and $i_0=i_1$. The length of a path $\gamma=(i_0,\dots,i_p)$ is the integer $p$. Note that there are no paths of length 0 and that every edge $(i,j)\in E$ is a path of length 1.

If $S\subset V$ is a subset of all the vertices, we will write $\overline S=V\setminus S$ and denote by $|S|$ the cardinality of the set $S$.

\begin{definition}
\emph{(Structural set).} Let $G=(V,E,w)\in\G$. A nonempty vertex set $S\subset V$ is a structural set of $G$ if 
\label{str-set}
\begin{itemize}
\item each cycle of $G$, that is not a loop, contains a vertex in $S$;
\item $w(i,i)\neq\lambda$ for each $i\in\overline S$.
\end{itemize}
\end{definition}

\begin{definition}
Given a structural set $S$, a \emph{branch} of $(G,S)$ is a path $\beta=(i_0,i_1,\dots,i_{p-1},i_p)$ such that  $i_0,i_p\in V$ and all $i_1,\dots,i_{p-1}\in\overline S$.
\end{definition}
We denote by $\mathcal{B}=\mathcal{B}_{G,S}$ the set of all branches of $(G,S)$.
Given vertices $i,j\in V$, we denote by $\mathcal{B}_{i,j}$ the set of all branches in $\mathcal{B}$ that start in $i$ and end in $j$. For each branch $\beta=(i_0,i_1,\dots,i_{p-1},i_p)$ we define the \emph{weight} of $\beta$ as follows:
\begin{equation}
w(\beta,\lambda):=w(i_0,i_1)\prod_{l=1}^{p-1}\frac{w(i_l,i_{l+1})}{\lambda-w(i_l,i_l)}.
\end{equation}

Given $i,j\in V$ set
\begin{equation}
R_{i,j}(G,S,\lambda):=\sum_{\beta\in\mathcal{B}_{i,j}}w(\beta,\lambda).
\end{equation}

\begin{definition}
\emph{(Isospectral reduction).} Given $G\in\G$ and a structural set $S$, the reduced adjacency matrix  $R_S(G,\lambda)$ is the $|S|\times |S|-$matrix with the entries $R_{i,j}(G,S,\lambda),i,j\in S$. This adjacency matrix $R_S(G,\lambda)$ on $S$ defines the reduced graph which is the isospectral reduction of the original graph $G$.
\end{definition}

\subsection{Isospectral Matrix Reductions}

For any matrix $M\in\W^{n\times n}$, let $N=\{1,2,\dots,n\}$. If the sets $R,C\subset N$ are nonempty, we denote by $M_{RC}$ the $|R|\times|C|$ submatrix of $M$ with rows indexed by $R$ and columns by $C$. Suppose that $S\subset N$ and its complement $\overline S=N\setminus S$ are nonempty. The isospectral reduction of $M$ over the set $S$ is defined as $$R_S(M,\lambda):=M_{SS}-M_{S\overline S}(M_{\overline S \overline S}-\lambda I)^{-1}M_{\overline SS}.$$
The only requirement for $S$ here is that the inverse matrix $(M_{\overline S\overline S}-\lambda I)^{-1}$ exists. This is a more general condition than that in Definition \ref{str-set}. Indeed for the isospectral graph reduction, there must be no non-loop cycles in $\overline S$, which means that after conjugating by a permutation $M_{\overline S\overline S}$ is a triangular matrix. Also, the weights of loops in $\overline S$ are not equal to $\lambda$. This ensures $M_{\overline S\overline S}-\lambda I$ is invertible, but this is a stronger condition. We would also like to point out that we cannot evaluate $R_S(M,\lambda)$ at points $\lambda_0\in\sigma(M_{\overline{SS}})$.

When both of these conditions hold the isospectral matrix reduction gives the same reduced matrix as the isospectral graph reduction (Theorem 2.1 \cite{webb14}). Isospectral reductions preserve the eigenvalues (Corollary 2.1 \cite{webb14}) and eigenvectors of a matrix (Theorem 1 \cite{torres14}, Theorem 3 \cite{shu18}). One can apply isospectral reductions sequentially and the final result only depends on the remaining indices/vertices (Corollary 2.3 \cite{webb14}).

\section{Spectral measurements}
\label{spec}

In this section we introduce several spectral measurements and facts for stochastic matrices.

\subsection{Definitions and properties}
A matrix $A\in\R^{n\times n}$ with entries $A=(a_{ij})$ is (column) \emph{stochastic} if $a_{ij}\geq 0$ and $\sum_{i=1}^n a_{ij}=1, \forall j=1,\ldots, n$. If $1$ is a simple eigenvalue of $A$ and there is no other eigenvalue of $A$ on the unit circle, we call $A$ \emph{non-critical}, otherwise $A$ is said to be \emph{critical}.
Let $\Delta^{n-1}$ be the simplex of all probability vectors
$$ \Delta^{n-1}:=\left\{ x=(x_1,\ldots, x_n)\in\R^n\colon x_i\geq 0, \; \sum_{i=1}^n x_i=1 \right\} .$$
For an integer $m\in\N$ we use the notation $a^m_{ij}$ for the entries of the power matrix $A^m$.
$A$ is \emph{primitive} if there exists $m\geq 1$ such that $a^m_{ij}>0, \forall i,j=1,\ldots, n$. Every primitive matrix is non-critical but the converse is not true.

The following concepts are taken from the theory of Markov Chains. See \cite{chung67} for details.

Given two vertices $i,j\in \{1,\ldots, n\}$,
we say that $i$ \emph{leads to} $j$ and write $i\rightsquigarrow j$ if there exists $m\geq 1$ such that $a^m_{ji}>0$. We say that $i$ and $j$ \emph{communicate} if $i\rightsquigarrow j$ and $j\rightsquigarrow i$, in which case we write $i \leftrightsquigarrow j$. A set of vertices $C\subseteq \{1,\ldots, n\}$ is called a class of $A$ if
\begin{enumerate}
    \item $i \leftrightsquigarrow j$ for all $i,j\in C$
    \item $C$ is saturated for $\leftrightsquigarrow$, i.e.,
    $i\in C$ and $i \leftrightsquigarrow j$ \, $\Rightarrow$ \, $j\in C$.
\end{enumerate}
 A class $C$ is \emph{essential} if for all $i\in C$ and $j\in\{1,\ldots, n\}$, $i\rightsquigarrow j$ implies $j\in C$.

A state $i\in \{1, \ldots, n\}$ is called \emph{recurrent} when for some $m\geq 1$, 
$a^m_{ii}>0$ and otherwise it is called a \emph{transient state}.

These concepts are known under different names in the Theory of oriented graphs, or digraphs. See Section 1.5 of \cite{jensen09}.
We present below a short dictionary of these terms for stochastic matrices $A$ with no transient states, where $G_A$ denotes the associated adjacency matrix.

\begin{center}
	\begin{tabular}{|c||c|}
		\hline
		Markov Chain Theory & Digraph Theory\\
		\hline \hline
		$A$ is irreducible & $G_A$ is strongly connected \\
		\hline
		$C$ is a class of $A$ & $C$ is a strong component of $G_A$ \\
		\hline
		$C$ is an essential  class of $A$ & $C$ is a terminal strong component of $G_A$ \\
		\hline		
	\end{tabular}
\end{center}
\bigskip

\begin{definition} \emph{(Inner spectral radius)}
  Let $\sigma(A)=\{\lambda_1,\lambda_2,\dots, \lambda_n\}$ be the eigenvalues of $A$ sorted by its absolute value in a way that $\lambda_1=1\ge|\lambda_2|\ge\dots \ge |\lambda_n|$. Then
 $$\rho_i(A):=\abs{\lambda_2}$$ is called the inner spectral radius of $A$.   
\end{definition}

\begin{definition}
For a stochastic matrix $A\in\R^{n\times n}$, a fixed point $q\in \Delta^{n-1}$ such that $A\, q=q$ is called a stationary measure of $A$.
\end{definition}

\begin{proposition}
\label{non-critical characterization}
For a stochastic matrix $A\in\R^{n\times n}$, the following are equivalent:
\begin{enumerate}
	\item $A$ is non-critical, i.e., $\rho_i(A)<1$,
	\item $A$ admits a unique essential class, which is aperiodic,
	\item There is a unique $v_\ast\in\Delta^{n-1}$ such that $A\, v_\ast=v_\ast$ and, moreover, $\lim_{m\to \infty} A^m\, v_0=v_\ast,\forall v_0\in\Delta^{n-1}$.
\end{enumerate}
\end{proposition}

\begin{proof}
This theorem follows from the classical theory of Markov Chains, 
see~\cite[Chapter V]{doob91}. We give a rough description of the results involved.

Given a stochastic matrix $A\in\R^{n\times n}$,
a vertex $i\in \{1,\ldots, n\}$ is transient if $i\not\leadsto i$. The relation $\leftrightsquigarrow$ is an equivalence relation on the set of non-transient vertices,
whose equivalence classes are precisely the classes defined above.
The set of all stationary measures is a compact polytope.
The extremal points of this compact convex set are the ergodic stationary measures,
where a stationary measure $q\in \Delta^{n-1}$ is said to be ergodic, respectively mixing, if the Markov shift determined by the pair $(A,q)$ is ergodic, respectively mixing, see~\cite{viana16}.
The support of a stationary measure is always a union of essential classes, i.e., stationary measures do not see transient vertices or non-essential classes. The map that assigns its support
to a stationary measure is a one-to-one correspondence between stationary measures and essential classes.

The partition of $\{1,\ldots, n\}$ in classes and transient vertices gives an upper triangular block representation of $A$.
Hence, the spectrum of $A$ is the union of the spectra of its essential classes and the spectrum of the submatrices corresponding to non-essential classes.
 
\end{proof}

\begin{corollary}
    For critical stochastic matrices, $\rho_i(A)=1$.
\end{corollary}

\bigskip

\begin{definition} \emph{(Diameter)}
For a stochastic matrix $A\in\R^{n\times n}$, 
$$\tau(A)  := \max_{i,j} \frac{1}{2}\,\sum_{k=1}^n\abs{ a_{ki}-a_{kj}}$$
is called the diameter of $A$.
\end{definition}

Consider the stochastic norm $\norm{x}_1:=\sum_{j=1}^n \abs{x_j}$.

Notice that $\tau(A)$ is half of the diameter of the image of the simplex $\Delta^{n-1}$ by   $A$ w.r.t. the norm $\norm{\cdot}_1$.
Moreover, the function $A\mapsto \tau(A)$ is a semi-norm on the space of matrices.

Since $$a_{ki}\wedge a_{kj}+\frac{1}{2}\, \abs{ a_{ki}-a_{kj}}=\frac{a_{ki}+a_{kj}}{2},$$
we get
$$\tau(A)=\max_{i,j}\sum_{k=1}^n[\frac{a_{ki}+a_{kj}}{2}-a_{ki}\wedge a_{kj}]
= 1- \min_{i,j} \sum_{k=1}^n a_{ki}\wedge a_{kj}.$$

\begin{proposition}
	For any stochastic matrix $A\in\R^{n\times n}$ and probability vectors $x,y\in \Delta^{n-1}$,
	$$ \norm{A\,x-A\,y}_1\leq \tau(A)\,\norm{x-y}_1 .$$
\end{proposition}

\begin{proof}
	Given two probability vectors $x,y\in \Delta^{n-1}$,
	averaging  we get
	\begin{align*}
		\norm{A\,x-A\,y}_1 &= \sum_{k=1}^n  \abs{ \sum_{i=1}^n a_{ki}\, x_i - \sum_{j=1}^n a_{kj}\, y_j  }\\
        &=\sum_{k=1}^n\abs{\sum_{i=1}^na_{ki}x_i\sum_{j=1}^ny_j-\sum_{j=1}^na_{kj}y_j\sum_{i=1}^nx_i}\\
		&\leq  \sum_{k=1}^n \sum_{i=1}^n  \sum_{j=1}^n \abs{  a_{ki} -  a_{kj} }\, x_i\, y_j \\
		&=  \sum_{i=1}^n  \sum_{j=1}^n \left( \sum_{k=1}^n \abs{  a_{ki} -  a_{kj} }\right) \, x_i\, y_j  \leq 2\, \tau(A).
	\end{align*}
	Next write $p=x-y$ with $p=p^+ -p^-$
	where $p^+:=(p^+_1,\ldots, p^+_n)$ and $p^-:=(p^-_1,\ldots, p^-_n)$ with
	$$ p^+_j:= \max\{x_j-y_j,0\} \; \text{ and } \;
	p^-_j:= \max\{y_j-x_j,0\} .$$
    Since $$\sum_{j=1}^n(x_j-y_j)=0=\sum_{j=0}^n(p^+_j-p^-_j)=\sum_{j=0}^np^+_j-\sum_{j=1}^np_j^-,$$ 
	we have $\alpha = \norm{p^+}_1=\norm{p^-}_1>0$. Applying the previous inequality to the probability vectors $\alpha^{-1}\,p^+$ and $\alpha^{-1}\,p^-$ 
	we get
	\begin{align*}
		\norm{A\,x- A\, y}_1 &=\norm{A\,(x-y)}_1 = \norm{A\,p^+ -A\, p^- }_1  \\
		&=\alpha\norm{A\alpha^{-1}p^+-A\alpha^{-1}p^-}_1\\
        &\leq 2\,\alpha\, \tau(A)   = \tau(A)\,( \norm{p^+}_1 + \norm{p^-}_1 )\\
		&=\tau(A)\,\norm{x-y}_1 .
	\end{align*}
\end{proof}

\begin{remark}
$\tau(A)$ is the operator norm of $A$ on
$\left( H, \norm{\cdot}_1\right)$, where  $H:=\left\{ x\in\R^n\colon \sum_{j=1}^n x_j=0 \right\}$.
In particular, $\rho_i(A)\leq \norm{A\vert_H}=\tau(A)$.
If $\tau(A)<1$ then $A$ is a contraction on $\Delta^{n-1}$.
\label{H}
\end{remark}

\begin{proposition}
    If $a_{ij}\geq c>0$ for all $i,j\in  \{1,\ldots, n\}$ then $$ \tau(A)\leq 1- n\, c <1 .$$
If there is an $i_0\in \{1,\ldots, n\}$ such that $a_{i_0 j}\geq c>0$ for all $j\in \{1,\ldots, n\}$
then $$ \tau(A)\leq 1- c <1 .$$
\end{proposition}

\begin{proof}
For instance, the second statement follows because $$ \tau(A) = 1- \min_{i,j} \sum_{k=1}^n a_{ki}\wedge a_{kj}\leq 1- \min_{i,j}  a_{i_0 i}\wedge a_{i_0 j}  \leq 1- c .$$
\end{proof}

\begin{definition}
    For any stochastic matrix $A\in\R^{n\times n}$, let $$m(A) = \min_{i, j}{a_{ij}}$$ be the smallest entry of $A$.
\end{definition}

Since $A$ is stochastic, if $m(A) = 1/n$, then all entries of $A$ are the same, exactly $1/n$. For a more general stochastic matrix whose entries are not all the same, we always have $m(A)<1/n$. 

If $A$ and $B$ are both matrices with non-negative entries, then we always have $m(A+B) \geq m(A) + m(B)$.

\begin{definition} \emph{(Spectral gap)}
    Let $g(A):=1-\rho_i(A)$ be the spectral gap for a stochastic matrix $A$. 
\end{definition}
Then $g(A)\ge1-\tau(A)$ and since $a_{ij}\ge m(A)$, we know $\tau(A)\le1-n\,m(A)$, therefore $g(A)\ge n\,m(A)$.

\subsection{Decrease of the semi-norm}
\begin{theorem}
Let $A\in \R^{n\times n}$ be a column stochastic matrix and $S\subset \{1,\ldots, n\}$. Denote by
$R_{S}(A,1)\in \R^{S\times S}$ the isospectral reduction of $A$ onto $S$ with value $1$ plugged in for $\lambda$. Then $R_S(A,1)$ is a stochastic matrix and
$$\tau(R_S(A,1))\le\tau(A).$$
\label{semi-norm}
\end{theorem}

\begin{proof}
We can prove this for the special case of $S = \{1, \ldots, n \} - \{r\}$ where we only remove one node in the reduction. Because sequential isospectral reductions are path-independent, we can always get to $R_S(A,1)$ for $|S|<n-1$ by removing one node at a time (Section 1.3 \cite{webb14}).

For $S=\{1,\dots,n\}-\{r\}$, let's write $R_S(A,1)$ as $R_r(A,1)$. This is the reduction where we remove node $r$ and plug in 1 for $\lambda$.Then $R_r(A,1)$ is a stochastic matrix. Indeed, the entries of $R_r(A,1)$ are

$$(R_r(A,1))_{ij} = a_{ij} + \frac{a_{ir}a_{rj}}{1-a_{rr}}\ge0, \quad i,j\neq r.$$

The column sum is 

$$\sum_{j\neq r}\left(a_{ij}+\frac{a_{ir}a_{rj}}{1-a_{rr}}\right)=\sum_{j\neq r}a_{ij}+a_{ir}=1.$$

Let $c_i(A)$ be column $i$ of $A$.

Then the semi-norm of $A$ is 
$$
 \tau(A) = \frac{1}{2} \sup_{i \neq j}{ \| c_i(A) - c_j(A) \|_1 }.
$$

Now let's estimate $\tau(R_r(A,1)$:
\begin{align*}
    \tau(R_r(A,1)) &= \frac{1}{2} \sup_{i \neq j} {\| c_i(R_r(A,1)) - c_j(R_r(A,1)) \|_1} \\
    &=\frac{1}{2} \sup_{i \neq j} {\sum_{l \neq r} \left| a_{li}+ a_{lr}(1-a_{rr})^{-1}a_{ri} -[a_{lj}+ a_{lr}(1-a_{rr})^{-1}a_{rj}]  \right|}\\
    &=\frac{1}{2} \sup_{i \neq j} {\sum_{l \neq r} \left| (a_{li} - a_{lj}) +  [a_{lr}(1-a_{rr})^{-1}a_{ri} -a_{lr}(1-a_{rr})^{-1}a_{rj}]  \right|}\\
    &=\frac{1}{2} \sup_{i \neq j} {\sum_{l \neq r} \left| (a_{li} - a_{lj}) +a_{lr}(1-a_{rr})^{-1} (a_{ri}- a_{rj})  \right|}\\
    &\leq\frac{1}{2} \sup_{i \neq j} {\sum_{l \neq r}  | (a_{li} - a_{lj})|  + \sum_{l \neq r}| a_{lr}(1-a_{rr})^{-1} ( a_{ri} - a_{rj})|}\\
    &=\frac{1}{2} \sup_{i \neq j} { \left\{\sum_{l \neq r} | (a_{li} - a_{lj})|+\left( \sum_{l \neq r}   a_{lr} \right) (1-a_{rr})^{-1} | a_{ri} - a_{rj}|\right\}}.
\end{align*}

Since $\sum_{l \neq r} a_{lr} = 1- a_{rr} $ the expression above becomes
$$
 \frac{1}{2} \sup_{i \neq j} { \left\{\sum_{l \neq r} | (a_{li} - a_{lj})|  +
 | a_{ri} - a_{rj}  |\right\}}  = \frac{1}{2} \sup_{i \neq j} { \left\{ \sum_{l=1}^n|(a_{li} - a_{lj})|  \right\}}
 = \tau(A).
$$
\end{proof}

\begin{lemma}
 Let $R_S(A,1)$ be the reduction of a stochastic matrix $A$ over set $S$ with value 1 for $\lambda$; then
$$
    m(R_S(A,1)) \geq \frac{m(A)}{1-|\overline{S}| m(A)}.
$$
\end{lemma}
\begin{proof}
We know that
\begin{align*}
   R_{S}(A,1) &= A_{S S} + A_{S \overline{S}} (I-A_{\overline{S} \overline{S}})^{-1} A_{\overline{S} S}\\
   &=A_{SS}+A_{S\overline{S}}A_{\overline{S} S}+A_{S\overline{S}}A_{\overline{S}\overline{S}}A_{\overline{S} S}+A_{S\overline{S}}A_{\overline{S}\overline{S}}^2A_{\overline{S}S}+\dots
\end{align*}
All the matrices in this expression are nonnegative.

Notice that for a given entry of $A_{S \overline{S}} A_{\overline{S} S}$ we have
$$
 (A_{S \overline{S}} A_{\overline{S} S})_{ij} \geq |\overline{S}| m(A)^2.
$$
For $A_{S \overline{S}} A_{\overline{S} \overline{S}} A_{\overline{S} S}$ we have
$$
  (A_{S \overline{S}} A_{\overline{S} \overline{S}} A_{\overline{S} S})_{ij} \geq |\overline{S}|^2 m(A)^3.
$$
and, more generally, $(A_{S \overline{S}} A_{\overline{S} \overline{S}}^k A_{\overline{S} S})_{ij} \geq |\overline{S}|^{k+1} m(A)^{k+2}$.
From those bounds we get
$$
 (R_S(A,1))_{ij} \geq m(A) + |\overline{S}| m(A)^2 + |\overline{S}|^2 m(A)^3 + \cdots = \frac{m(A)}{1-|\overline{S}| m(A)}.
$$
Here we are assuming that $|\overline{S}| m(A) < 1$; but this is true since $|\overline{S}| < n$ and $m(A) \leq 1/n$.
\end{proof}
When $m(A)$ is positive we get that $m(R_S(A,1)) > m(A)$. Moreover when $m(A)$ is large, the image of the cone of non-negative vectors under the stochastic matrix $A$ is smaller, or $A$ is more contractive. Since $m(R_S(A,1))$ is larger than $m(A)$ we get that $R_S(A,1)$ is indeed more contractive than $A$.

\subsection{Gershgorin estimate}

For a matrix $A\in\C^{n\times n}$ we define for each $i \in \{1, \ldots, n\}$ the $i-$th absolute row sum $$r_i(A) = \sum_{j \neq i} |a_{ij}|.$$ 

\begin{theorem}[Gershgorin \cite{gershgorin31}]
Let $A\in\C^{n\times n}$. Then all eigenvalues of $A$ are contained in the set
$$\Gamma(A)=\bigcup_{i=1}^n\{\lambda\in\C:|\lambda-a_{ii}|\le r_i(A)\}.$$
\end{theorem}

Now let us assume that $A$ has positive elements only and is both column and row stochastic, i.e. all its row and column sums are 1.

We define the matrix $R_n(A,1)$ as the isospectral reduction where vertex $n$ is eliminated. Then each entry of the reduction is given by
$$
   R_n(A,1)_{ij} = a_{ij} + \frac{a_{in} a_{nj}}{1 - a_{nn} }
$$
for $i, j \in \{1, \ldots, n-1\}$. Each entry of $R_n(A,1)$ is still positive and one can check that the row sums as well as column sums are still $1$.

Let us compute the first absolute row sum $r_1(R_n(A,1))$.
\begin{align*}
 r_1(R_n(A,1))=& R_n(A,1)_{12} + \ldots + R_n(A,1)_{1, n-1} =\\
  &a_{12} +\ldots + a_{1, n-1} + \frac{a_{1n}}{1-a_{nn}} \left( a_{n2} + \ldots + a_{n, n-1}\right).
\end{align*}
Notice that $ a_{n2} + \ldots + a_{n, n-1}  =  1- a_{n1} - a_{nn}$; hence the expression above becomes
\begin{align*}
 r_1(R_n(A,1))=&a_{12} + \ldots + a_{1, n-1} + \frac{a_{1n}}{1-a_{nn}} \left(  1- a_{n1} - a_{nn}   \right) \\
 <&a_{12} + \ldots + a_{1, n-1} +  a_{1n} = r_1(A).
\end{align*}

Hence $r_i(R_n(A,1)) < r_i(A)$  when $i=1$ and we can show that this is true for $i= 2, 3, \ldots, n-1$ as well. This agrees with the results in \cite{webb12}.

\section{Computation methods and a new approach}
\label{comp}

In this section we will give a brief explanation of Gaussian elimination, Perron-Frobenius iterations and propose a new approach to compute the stationary measure of a large Markov chain, based on isospectral reductions.

\subsection{Iterative methods}

Let $A$ be some $n\times n$ stochastic matrix.
In general, if $A$  is  non-critical  we can approximate its stationary measure by the iterates $A^m\, v_0$ starting from any vector $v_0\in\Delta^{n-1}$ (Proposition~\ref{non-critical characterization}). The Perron-Frobenius method returns  an approximation of the stationary measure of $A$ with an error up to $10^{-p}$, for some positive integer $p$.

\begin{pseudocode}[ruled]{PerronFrobenius}{A, p}
	\label{PerronFrobenius}
    \LOCAL{v_0, v_1, \Delta}\\
		
	v_1:= \text{ random point in }  \Delta^{n-1}\\	
	
	\REPEAT  
		v_0:=v_1 \\
		v_1:= A\, v_0\\
		\Delta:= v_1-v_0 \\

	\UNTIL
		\norm{\Delta}^2 <10^{-2 p} \\
		
	\RETURN   { v_1  }
\end{pseudocode}

The iterates $v_m:=A^m\, v_0$ converge to the fixed point
$v_\ast=A\,v_\ast$ at a geometric rate dependent on the inner spectral radius of $A$,
$$\abs{v_m-v_\ast } \lesssim [\rho_i(A)]^m.$$ So the number of iterations needed to reach precision $10^{-p}$ is $-p\ln10/\ln\rho_i$ while each iteration takes $n^2$ operations. In total the computational cost is $-n^2p\ln10/\ln\rho_i$.
 Notice that in general, when $\rho_i(A)\ll 1 -\frac{1}{n}$, the computational cost can be of order $O(n^2)$, much better than Gaussian elimination, which is of order $O(n^3)$. On the other hand if $\rho_i(A)\gtrsim 1-\frac{1}{n}$ the cost rises up to at least $O(n^3)$, the same order as in Gaussian elimination. In fact, in this case  Gaussian elimination will likely generate ill conditioned matrices leading to  inaccurate answers.

The most used methods to approximate an eigenvector (of a general non normal matrix) are iterative methods and among them variations of the Arnoldi method~\cite{saad11}.
All these methods for approximating the eigenvector associated with a given eigenvalue work poorly in the presence of a small spectral gap measured from the given eigenvalue~\cite{eisenstat98}.

\subsection{Isospectral algorithmic scheme}

When a stochastic matrix $A$ has more than one eigenvalue very close to $1$, no standard iterative method will work well to approximate the stationary measure of $A$.   This is precisely the case when  Isospectral Theory can be useful.

The isospectral reduction \begin{equation}
\label{Reduction equation}
   R_S(A,1)=A_{SS}-A_{S\overline{S}}(A_{\overline{SS}}-I)^{-1}A_{\overline{S}S} 
\end{equation} is a stochastic matrix if $A$ is stochastic (Theorem \ref{semi-norm}). In fact, the stationary vector $v_R$ of $R_S(A,1)$  is the projection of the stationary measure of $A$ to the coordinates in the set $S$.  We can also reconstruct $v_A$ from $v_R$ as follows (\cite{torres14,shu18}).

\begin{align}
    v_A=\begin{bmatrix}
        v_R\\
        -(A_{\overline{SS}}-I)^{-1}A_{\overline{S}S}v_R
    \end{bmatrix}
\end{align}

 This suggests an alternative scheme to approximate the stationary measure of $A$. 

\begin{pseudocode}[ruled]{Isospectral}{A, p}
	\label{Isospectral}	
	\LOCAL{S, v_R, v_A}\\
    \text{(i) Choose a subset } S \text{ of } \{1,\ldots, n\} \\	
    \text{(ii) Reduce } A \text{ over } S \text{ to get a smaller stochastic matrix } R_S(A,1)\\
    \text{(iii) Compute the stationary measure } v_R \text{ of } R_S(A,1), \text{ with an error up to } 10^{-p} \\
    \text{(iv) Reconstruct the stationary measure } v_A \text{ of } A \text{ from } v_R\\
	\RETURN   { v_A  }
\end{pseudocode}

There are several options for steps (i), (ii) and (iii). 

\begin{enumerate}
    \item Regarding the first step, besides choosing the size $s=|S|$, we can:
\begin{itemize}
    \item  choose $S=\{1,\dots,s\}$,
    \item  choose $S\subset \{1,\ldots, n\}$ randomly of size $s$,
    \item choose $S\subset \{1,\ldots, n\}$ in some deterministic way that attempts to minimize the cost of step (iii).
\end{itemize}

\item Regarding the second  step there are two possibilities:
\begin{itemize}
    \item Perform the isospectral reduction in one-step using equation ~\eqref{Reduction equation}. The down side of this approach is the possibility of an ill conditioned inversion 
    $(I-A_{\overline{SS}})^{-1}$.
    Based on empirical evidence, randomizing the choice of $S$ seems a good way to address this issue. The reduction~\eqref{Reduction equation} takes 
    $$O\left((n-s)^3+(n-s)^2s+s^2(n-s)+s^2\right)$$
    multiplications\footnote{because they are much less expensive we don't count additions.} to compute, where the first term corresponds to  inversion of the matrix $I-A_{\overline{SS}}$ and the others to the multiplications of the matrices involved. Also choosing $s$ random nodes takes $O(s)$ time.

\bigskip

    \item Perform the isospectral reduction in $n-s$ steps where at each step the isospectral reduction is performed by  choosing and removing a single node. For instance in the first step,  choosing the node $k$, so that 
$S=\{1,\ldots, n\}\setminus\{k\}$,
the reduced matrix has entries
$$ r_{ij}=a_{ij}+\frac{a_{ik}\, a_{kj}}{1-a_{kk}} \qquad i,j\in S.$$
Avoiding choices of nodes where $a_{kk}\approx 1$ the problem of computing the isospectral reduction is always well conditioned. By~\cite[Theorem 2.5]{webb14}, the order in
$\{1,\ldots, n\}\setminus S=\{k_1,\ldots, k_{n-s}\}$ by which the one-step reductions are performed is indifferent in the sense that they all lead to the same matrix $R_S(A,1)$ in~\eqref{Reduction equation}.
The \emph{time cost} of the step by step isospectral reduction is 

\begin{align*}
    \sum_{j=s+1}^{n}   & (j-1)^2+j-1 =\frac{(n+1) n (n-1) - (s+1) s (s-1) }{3}\\
    & = O\left( \;  \frac{n^3-s^3}{3}\; \right) =O\left( \; \frac{(n-s)^3}{3} + (n-s)^2 s +(n-s) s^2+s \; \right).
\end{align*}
\end{itemize}

Notice that this \textit{time cost} is of the same order as the single step reduction time cost, and even better unless  the computational cost of $(I-A_{\overline{SS}})^{-1}$ gets to be less than $(n-s)^3/3$. The above considerations show that the single isospectral reduction from an $n\times n$ matrix to one of size $s\times s$ in temporal cost terms is computationally equivalent to performing a sequence of $s$ single one step reductions. In this comparison we are not taking into account spatial memory costs.

\bigskip

\item Finally, regarding the third step  there are also several options:
\begin{itemize}
    \item  Use one of the available iterative methods to approximate the dominant eigenvector of the reduced matrix $R_S(A,1)$.
    \item  After reduction, estimate the spectral gap $\rho_i(R)$ and if it is still very close to $1$ repeat steps (i) and (ii). Otherwise
    apply one of the existing iterative methods to approximate the dominant eigenvector of the reduced matrix $R_S(A,1)$.
\end{itemize}

\end{enumerate}

In the next subsection we show results of numerical experiments that illustrates the advantage of combining traditonal iterative methods with isospectral reduction through the proposed algorithmic scheme.


\subsection{Comparison of the methods}



\newcommand{\Prob}{\mathrm{Prob}}
\newcommand{\SM}{\mathcal{M}}

Let $P:\R^n_+\to \Delta^{n-1}$
be the projection
$$ P(x_1,\ldots, x_n):=  \left(\frac{x_1}{\sum_{k=1}^n x_k}, \ldots, \frac{x_n}{\sum_{k=1}^n x_k} \right),  $$
which extends to the  projection $P:\R_+^{n\times n}\to \SM_n$ onto the space of all $n$ by $n$ stochastic matrices, where $P(A)$ has entries
$$ P_{ij}(A):= \frac{a_{ij}}{\sum_{k=1}^n a_{kj}}. $$

Given a probability measure $\mu\in\Prob(\R_+)$, let $\mu^n\in\Prob(\R_+^n)$ be the Cartesian product measure
and $\Pp_\mu:= P_\ast \mu^n\in\Prob(\SM_n)$ be the push forward probability on the space of stochastic matrices.

The probability $\mu$ is called \textit{symmetric} if it is invariant
under the involution $x\mapsto \frac{1}{x}$, or if
$$\mu(0, x)=\mu(\frac{1}{x},\infty),\forall x\in (0,\infty).$$

We say that $\mu$ has a \textit{heavy tail} if there exists
$\alpha\in (0,1)$ such that 
$$\lim_{x\to \infty} x^\alpha\, \mu(x,\infty)=c>0 .$$

\begin{example}[Burr distribution]
The Burr distribution 
is a family of symmetric laws $\mu_\alpha\in\Prob(\R_+)$  with heavy tail,  cumulative distribution function,
$$ F_\alpha(x):=1-\frac{1}{1+x^\alpha} $$ 
	and probability density function,
$$ f_\alpha(x):=\frac{\alpha}{x^{1-\alpha}\, (1+x^\alpha)^2} . $$
\end{example}

We will use Burr distribution in generating some of the stochastic matrices in the following numerical experiments. A stochastic matrix $A$ generated with Burr distribution is conjectured to have its spectral radius $\rho_i(A)\sim\sqrt{1-\alpha}$ \cite[Remark 1.3]{bordenave17}.

\begin{experiment}[Nested reductions - Burr distribution 1]

In this experiment we randomly generate a $1000\times1000$ stochastic matrix $A$ and compare the default iterative method in Mathematica with random Isospectral method (Algorithm \ref{Isospectral}). The computation error is $E=\|A\mathbf{v}-\mathbf{v}\|$, where $\mathbf{v}$ is the computed stationary measure.

\begin{figure}[h]
    \centering
    \includegraphics[scale=0.7]{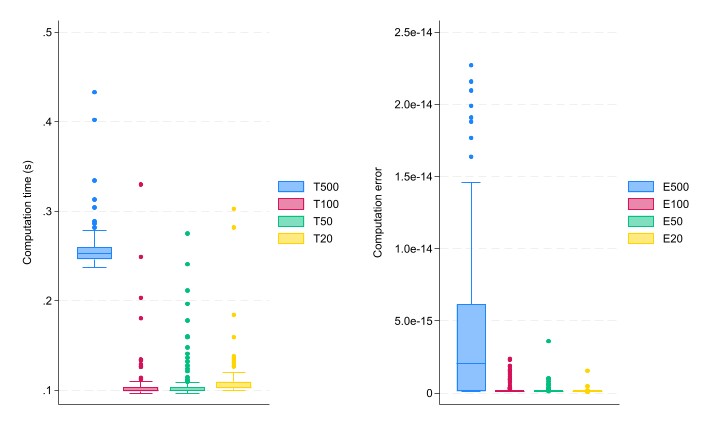}
    \caption{Computation times and errors of nested random isospectral reduction and reconstruction of the stationary measure. In this specific example, the Mathematical default method does not converge.}
    \label{diverge}
\end{figure}

The matrix $A$ was chosen to have $100$ nonzero entries in each column, which follow an i.i.d. Burr distribution with $\alpha=0.1$. We computed $\rho_i(A)=1$ (by Mathematica default accuracy) and the Mathematica default method failed to converge, stopping the computation after $T=5.17114$ seconds with error $E=1.43583\times10^7$. Then we used Algorithm \ref{Isospectral}, with one step reduction and the Mathematica default method to find the dominant eigenvector of the reduced matrix,  to compute the stationary measure with structural sets $S_1\supset S_2\supset S_3\supset S_4$, randomly chosen with sizes $|S_1|=500,|S_2|=100,|S_3|=50,|S_4|=20$. We recorded the corresponding computation times $T500, T100, T50, T20$ and computation errors $E500,E100,E50,E20$ for 100 different random subset sequences $S_1, S_2, S_3, S_4$ in Figure \ref{diverge}. As one can see, all computation times are within half a second and accuracies are at the order of $10^{-14}$ or higher. The errors are quite small. Among all these fairly fast and accurate computations, the size $500$ reduction performed the poorest and yet it's still a great improvement over the default Mathematica method.
\end{experiment}

\begin{experiment}[Nested reductions - Burr distribution 2]
In the second experiment we generated a different $A$, still a $1000\times1000$ stochastic matrix, now with 80 nonzero entries in each column with i.i.d. Burr distribution ($\alpha=0.5$). The new inner spectral radius is $\rho_i(A)=0.999808$. Mathematica took 0.415732 seconds to find the stationary measure of $A$ and the error was $1.97685\times10^{-8}$. We repeated the same procedure as in the previous experiment, reducing $A$ to randomly chosen nested sequences of subsets of sizes 500, 100, 50 and 20 and recording the computation times and errors, see Figure \ref{burr3}. Again, the size $500$ reduction performs the poorest. The computation time is comparable to Mathematica but the accuracy is still better. The other reductions perform even better, faster and more accurate.

\begin{figure}[h]
    \centering
    \includegraphics[scale=0.7]{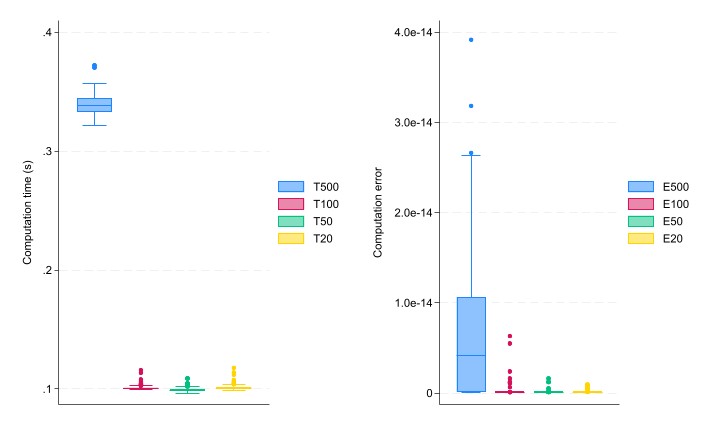}
    \caption{Computation times and errors of nested random isospectral reduction and reconstruction of stationary measure. The results are consistently more accurate.}
    \label{burr3}
\end{figure}

Just to give more details about the results, we also show a scatter plot of $E50$ vs $T50$ in Figure \ref{independence-3}. For this specific stochastic matrix $A$, and reducing to 100 random subsets of size 50, it seems the accuracy is fixed while the computation time varies.

\begin{figure}[h]
    \centering
    \includegraphics[scale=0.7]{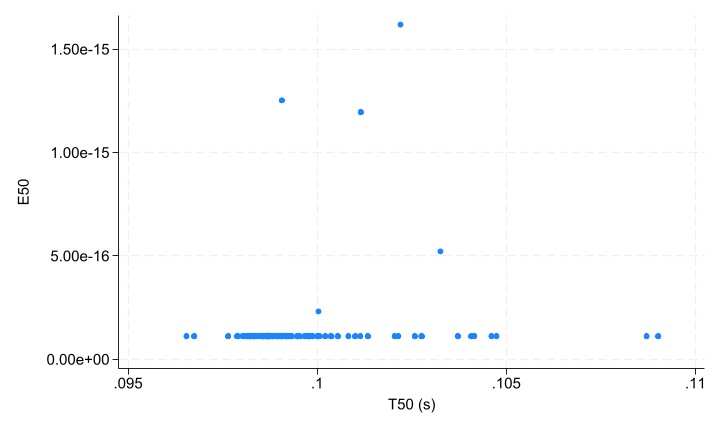}
    \caption{Scatter plot of $E50$ vs $T50$ for the second sparse Burr stochastic matrix.}
    \label{independence-3}
\end{figure}
\end{experiment}

One might have noticed in these two examples the inner spectral radius of the original matrix is either 1 or very close to 1. This is exactly where the Mathematica default method struggles. One might wonder how the reductions perform when the inner spectral radius is far from 1.

\begin{experiment}[Nested reductions - Uniform distribution]
We generated a $1000\times1000$ stochastic matrix $A$ with 100 nonzero entries in each column that follow i.i.d. Uniform distribution. Now $\rho_i(A)=0.111986$, and the default Mathematica method took 0.171988 seconds to find the stationary measure with an error of $1.28664\times10^{-17}$, which is great performance. We again reduced $A$ to random sequences of subsets of sizes $500,100,50$ and $20$, see Figure \ref{uniform}. Now the size 500 reductions are slower than the Mathematica default and have comparable accuracy to the default method. The rest of the reductions are slightly faster and slightly more accurate.

\begin{figure}[h]
    \centering
    \includegraphics[scale=0.7]{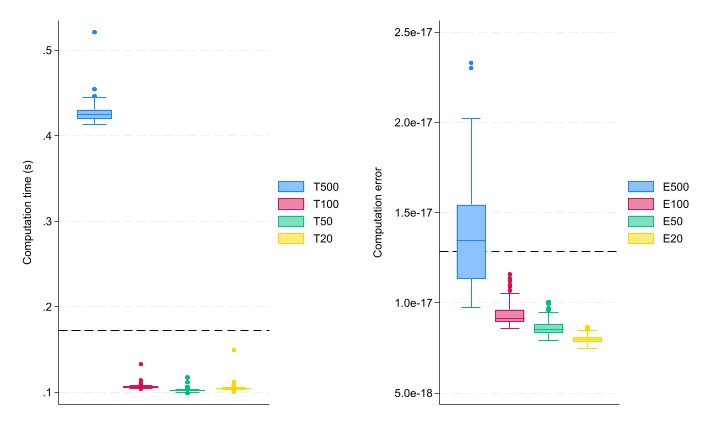}
    \caption{Computation time and error of nested random isospectral reduction and reconstruction of stationary measure. In this specific example the nonzero entries of the stochastic matrix are generated with uniform distribution. The dashed lines are the performance of the Mathematica default method.}
    \label{uniform}
\end{figure}
\end{experiment}

While computing the reduction and reconstruction takes time, for a stochastic matrix $A$ whose inner spectral radius is very close to 1, typically $\rho_i(R_S(A,1))$ is less than $\rho_i(A)$. This makes the computation of the stationary measure of $R_S(A,1)$ faster than that of $A$. Obviously $R_S(A,1)$ is also a smaller stochastic matrix compared to $A$. Even if $\rho_i(R_S(A,1))=\rho_i(A)$, computing the stationary measure of $R_S(A,1)$ is still faster than that of $A$. This explains why when the subset sizes are small ($\lesssim10\%$ of the original size), even for uniform distribution generated stochastic matrices (whose inner spectral radii are already small), we still have some gain in both computation time and accuracy.

\begin{experiment}[Barabási–Albert graphs - fluctuations]

In this experiment we focused on Barabási–Albert graphs \cite[Chapter 5]{cohen10}. First we generate a random digraph with the Barabási–Albert model, with 1000 vertices where a new vertex with degree $3$ was added at each step of the construction. Edges were rewired with probability 0.1. This graph was produced with  the {\rm IGBarabasiAlbertGame} and {\rm IGRewireEdges} functions from the
{\em IGraph} interface for Mathematica, see~\cite{IGraph2024}.
Edges of the resulting digraph were assigned i.i.d. random weights with Burr distribution ($\alpha=0.5$). Normalizing the columns of the weighted adjacency matrix gave  us a stochastic matrix $A$.

In this example we found $\rho_i(A)=0.99958$. The Mathematica default method took 0.212462 seconds to find the stationary measure of $A$ with an error of size $3.61832\times10^{-9}$. Here we reduced the matrix $A$ one-hundred times over random subsets $S\subset\{1,2,\dots,1000\}$ of size 100 and recorded the computation time and error of the Isospectral method. See Figure \ref{table4}. We can see the computation times concentrate around 0.085 seconds, and the errors spread out more than the computation times. But overall our method is faster and more accurate than the Mathematica default, since again the inner spectral radius is very close to 1.

\begin{figure}[h]
    \centering
    \includegraphics[scale=0.7]{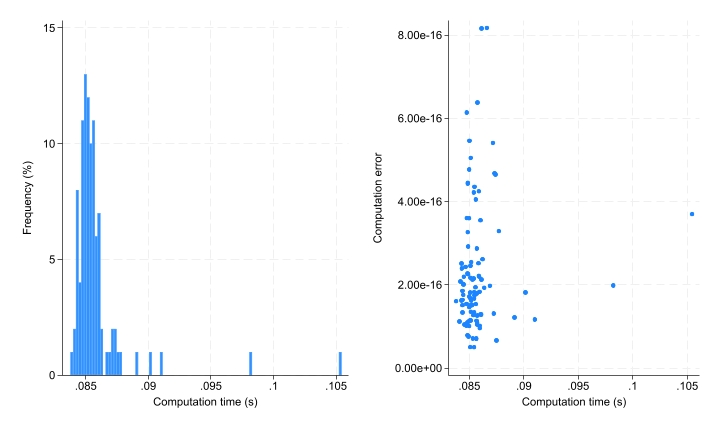}
    \caption{Fluctuation of computation times and errors for our method with reduction and reconstruction using isospectral reductions applied to a Barabási–Albert graph.}
    \label{table4}
\end{figure}
\end{experiment}
\begin{experiment}

Now that we have an idea of the fluctuation in performance of our method on a fixed stochastic matrix associated with a Barabási–Albert graph, we are ready to vary the matrix $A$. For 100 times, we generate a Barabási–Albert graph with 100 vertices, where a new vertex with degree 3 is added at each step of the construction. No edges are rewired. Then we assign i.i.d. weights with Burr distribution ($\alpha=0.5$) to edges and normalize to get a stochastic matrix $A$. We compare the performance of the Mathematica default method with our method reducing $A$ over a random subset of size 10. Let $T_M$ be the computation time of the Mathematica default method and $T_I$ be the computation time of the method using Isospectral reductions. Similarly $E_M$ and $E_I$ are computation errors of the corresponding methods. We show the scatter plot of $\log_{10}(E_M/E_I)$ against $T_I/T_M$. We repeat this procedure with all the same parameters except the graph size increases to 500 and the subset size to 50. See Figure \ref{BA100-500-1000}. Our errors are consistently $10^{-6}$ times below the errors of the Mathematica default. For the most part our computation is faster as well, with some outliers where the computation time is comparable to the Mathematica default method.

\begin{figure}[h]
    \centering
    \includegraphics[scale=0.7]{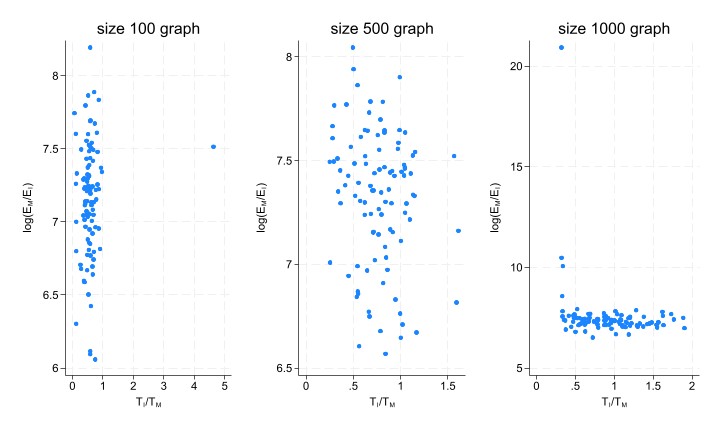}
    \caption{Computation performance comparison between the Mathematical default method and Isospectral reduction scheme on Barabási–Albert graphs.}
    \label{BA100-500-1000}
\end{figure}

We run this procedure for size 1000 graphs as well with all the same parameters and reducing to size 100 subsets. The result is in Figure \ref{BA100-500-1000}. We can see the errors of our method are $10^{-6}$ times below the errors of the default method and the computation time is comparable to the default method. If one wants to reach the same level of accuracy with the Mathematica default method, they would take longer time, not necessarily faster than our method.
\end{experiment}
\begin{experiment}

Lastly, we generate a size 1000 Barabási–Albert graph where a new vertex with degree 3 is added at each step of the construction. Edges are rewired with probability 0.1. Again we assign i.i.d. weights with Burr distribution ($\alpha=0.5$) to edges and normalize the columns to get a stochastic matrix $A$. Here we have $\rho_i(A)=0.998474$. We choose 500 random subsets $S$ of size 50 and reduce $A$ to $R_S(A,1)$ and record the ratios $\rho_i(R_S(A,1))/\rho_i(A)$ (Figure \ref{8-10}). We can see $\rho_i(R_S(A,1))/\rho_i(A)\le1$ but it's quite concentrated near 1. The average ratio is around 0.914696 with the 75 percentile at 0.98186. We would like to point out that these matrices generated with Barabási–Albert model and Burr distribution are on average very ill conditioned and we frequently see warnings from Mathematica.

\begin{figure}[h]
    \centering
    \includegraphics[scale=0.7]{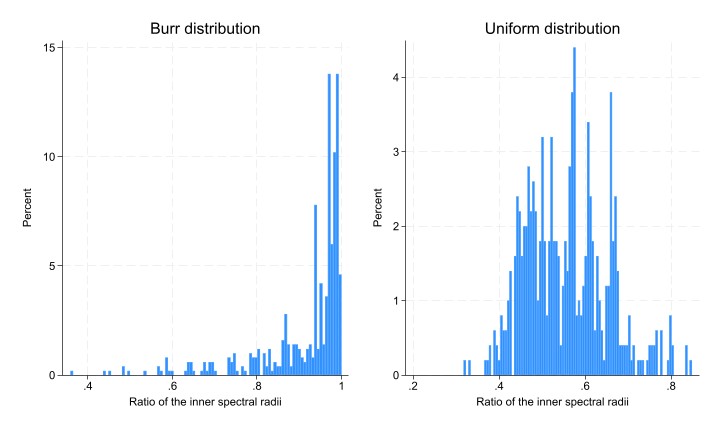}
    \caption{Computation result for size 1000 Barabási–Albert graphs and the ratio of the inner spectral radii before and after reduction.}
    \label{8-10}
\end{figure}

To see the effects of Barabási–Albert graphs alone we repeat this experiment changing Burr distribution to Uniform distribution. The corresponding inner spectral radius is $\rho_i(A)=0.794564$. From Figure \ref{8-10} we can see now the ratio of the inner spectral radii is more spread out compared to Burr distribution and on average concentrated around 0.6. This shows the inner spectral radius does generally decrease after isospectral reductions, with a high probability.
\end{experiment}

In our experiments which do not involve the Barabási–Albert model, edges are generated with a uniform distribution. When the weights of these edges follow Burr distribution, our method performs very well. In the experiments with Barabási–Albert models, the edges already follow a power law or non-uniform distribution, it seems now if we additionally make the weights of these edges follow Burr distribution, the advantage of our method is not as obvious as when the edge weights follow uniform distribution.

\section{Matrices that become positive after reduction}
\label{sec-positivematrices} 

Intuitively, we know that if all the entries of a stochastic matrix are close to each other, its semi-norm is small. If every entry is the same, we have a semi-norm of 0 and can find the stationary measure after one step of iteration. Meanwhile for a stochastic matrix where two columns have interlacing zero entries we have a semi-norm of 1 and probably need many iterations to approximate the stationary measure. 

We would like to see fewer zero entries in the matrix after reduction, ideally no zero entries, and we have a strictly positive matrix whose semi-norm is small. We prove some general results and construct specific examples where the matrix becomes strictly positive after reduction in the following subsections.

\subsection{General results}

\begin{proposition}
 Let $m\in\{1,2,\dots,n\},n\ge3$, $S=\{1,2,\dots,n\}\setminus\{m\}$. $A$ is an  $n \times n$ non-negative irreducible matrix and $\lambda_0\ge0$ is the dominant eigenvalue of $A$, i.e. $\lambda_0$ is an eigenvalue of $A$ and any other eigenvalue $\lambda$ (possibly complex) in absolute value is not larger than $\lambda_0$, $|\lambda|\le\lambda_0$. If
 \begin{enumerate}
     \item each row and each column of $A$ has at most one element that is zero,
     \item column $m$ has the lowest column sum,
     \item and row $m$ is positive, that is, it only has positive entries,
 \end{enumerate}
 then $R_S(A,\lambda_0)$ is a positive matrix. This matrix is also called the isoradial reduction of $A$ \cite{webb20}.
\end{proposition}
\begin{proof}
    Matrix $A$ has an eigenvector with non-negative components corresponding to the dominant eigenvalue $\lambda_0$, which is larger than the absolute value of any other eigenvalue.  Eigenvalue $\lambda_0$ satisfies

    $$\min_{i}{ \left(\sum_{j=1}^n a_{ij}   \right) }  \leq \lambda_0 \leq\max_{i}{ \left(  \sum_{j=1}^n a_{ij}\right)}.$$

    Since $m \in \{1, 2, \ldots, n \}$ is where the minimum above is reached and $n\ge3$,

    $$a_{mm} <  \sum_{j=1}^n a_{mj} \leq \lambda_0   \;\; \Rightarrow \lambda_0 - a_{mm} > 0.$$
    The entries of $R_S(A,\lambda_0)$ are
    $$(R_S(A,\lambda_0))_{ij}  = a_{ij} +  \frac{a_{im}a_{mj}}{\lambda_0-a_{mm}},\quad \forall i,j\in S. $$
    Since $a_{mj}>0$ for all $j$, the term $\dfrac{a_{mj}}{\lambda_0-a_{mm}}$ is always positive.  If
    $a_{im}$ is zero then $a_{ij}$ (also on row $i$) has to be positive. If $a_{ij}$ is zero then $a_{im}$ is positive.  Either way we have $ (R_S(A,\lambda_0))_{ij} > 0,\forall i, j\in S$.
\end{proof}

Notice that if $S=\{i\}$, a single node, the reduction $R_S(A,\lambda_0)=\lambda_0$ is positive. If we remove one node at a time, we will eventually reach a positive matrix. The above proposition could be used to find large structural sets $S$ such that the reduction is positive.

\begin{proposition}
    Let $S\subset\{1,\dots,n\}$ and let $A$ be an $n\times n$ stochastic matrix. If 
    \begin{enumerate}
        \item $a_{ii}<1,\forall i=1,\dots,n$,
        \item each row and each column has at most $m$ elements that are zero,
        \item and the cardinality of the set of removed vertices $|\overline{S}|>2(m-1)$,
    \end{enumerate}
    then $R_S(A,1)$ is positive.
    \begin{proof}
        The entries of $R_S(A,1)$ are
        $$
        (R_S(A,1))_{ij} = a_{ij} + \sum_{k \in \overline{S}} a_{ik} (1-a_{kk})^{-1} a_{kj}  + \dots, \forall  i, j\in S.
        $$
        Here the omitted terms are higher order terms like $a_{ik}(1-a_{kk})^{-1}a_{kl}(1-a_{ll})^{-1} a_{lj}$ and they are non-negative.
        
        If $a_{ij}$ is positive we are done.
        
        If $a_{ij}$ is zero, then at most $m-1$ of the $a_{ik}$ ($k \in \overline{S}$) terms can  be
        zero. Similarly, at most $m-1$ of the $a_{kj}$ ($k \in \overline{S}$) terms can be zeros. Since $|\overline{S}|>2(m-1)$ then
        we have at least one non-zero element in the sum above, showing that $(R_S(A,1))_{ij} > 0$.
    \end{proof}
\end{proposition}

\begin{remark}
    If $A\in\R^{n\times n}$ is a non-negative matrix with at most $m$ zero entries on each row and each column, and $n>2m$, then $A^2$ is positive.
\end{remark}

Indeed,
$$
  (A^2)_{ij} = \sum_{k=1}^n a_{ik} a_{kj}
$$
and we have at most $2m < n$ vanishing terms in the sum above. Every entry of $A^2$ is positive.

\begin{proposition}
    Consider the non-negative matrix
    $$A =\begin{bmatrix}A_{\overline{S}\overline{S}}  &  A_{\overline{S} S}  \\
     A_{S \overline{S}} &  A_{S S}\end{bmatrix}.$$
     Here $A$ is column stochastic, $A_{\overline{S}\overline{S}}$ is primitive, i.e., there exists some integer $l$ such that $(A_{\overline{S}\overline{S}}^l)_{ij} > 0$ for all $i, j \in \overline{S}$. If in addition, 
    \begin{enumerate}[(a).]
        \item the sum of each column of $A_{\overline{SS}}$ is strictly less than $1$.
        \item  for $i, j \in S$   either $(A)_{ij} > 0$ or the column $j$ of
        $A_{\overline{S} S} $ has at least one positive element;
        \item  for $i, j \in S$   either $(A)_{ij} > 0$ or the row $i$ of
        $A_{S \overline{S}}$ has at least one positive element;
    \end{enumerate}
    then $R_S(A,1)$ is positive.
\end{proposition}
\begin{proof}
    Observe that
    \begin{align*}
   R_{S}(A,1) &= A_{S S} + A_{S \overline{S}} (I-A_{\overline{S} \overline{S}})^{-1} A_{\overline{S} S}\\
   &=A_{SS}+A_{S\overline{S}}A_{\overline{S} S}+A_{S\overline{S}}A_{\overline{S}\overline{S}}A_{\overline{S} S}+A_{S\overline{S}}A_{\overline{S}\overline{S}}^2A_{\overline{S}S}+\dots
\end{align*}
This converges because of (a). Since $A_{\overline{SS}}$ is primitive, $A_{\overline{SS}}^l$ is positive for some $l$. (b) and (c) ensure either $(A_{SS})_{ij}>0$ or $(A_{S\overline{S}}A_{\overline{SS}}^lA_{\overline{S}S})_{ij}>0$.
\end{proof}

Under the conditions of this proposition the reduction of $A$ over $\overline{S}$
$$R_{\overline{S}}(A,1) = A_{\overline{S}\overline{S}} + A_{\overline{S} S} (I-A_{S S})^{-1} A_{S\overline{S}}$$
is primitive when it exists.

In the final discussion of this subsection we consider matrices whose entries are close to $1/n$. Notice that the matrix whose entries are exactly all $1/n$ strongly contracts the hyperplane $H$ in Remark \ref{H}. We try to estimate how far the reduced matrix is from the averaging matrix whose entries are all identical.

Assume the entries of the stochastic matrix $A\in\R^{n\times n}$ satisfy
$$
a_{ij} \leq d(n) = \frac{1}{n} + \epsilon(n).
$$

The reduction of $A$ with eigenvalue 1 is
$$
R_S(A,1) = A_{SS} + A_{S \overline{S}} A_{\overline{S} S} + A_{S \overline{S}} A_{\overline{S} \overline{S}} A_{\overline{S} S} + A_{S \overline{S}} A_{\overline{S} \overline{S}}^2 A_{\overline{S} S} + \cdots
$$

The entries of the terms in the sum satisfy
\begin{align*}
    &(A_{S \overline{S}} A_{\overline{S} S})_{ij} \leq |\overline{S}| d(n) d(n),\\
    &(A_{S \overline{S}} A_{\overline{S} \overline{S}} A_{\overline{S} S})_{ij} \leq |\overline{S}| |\overline{S}| d(n) d(n) d(n),\\
    & \dots \\
    &(A_{S \overline{S}} A_{\overline{S} \overline{S}}^m A_{\overline{S} S})_{ij} \leq |\overline{S}|^{1+m}  d(n)^{2+m}.
\end{align*}
Hence the entries of the reduction satisfy
$$
  (R_S(A,1))_{ij} \leq d(n) + |\overline{S}| d(n)^2 + |\overline{S}|^2 d(n)^3 + \ldots = \frac{d(n)}{1-|\overline{S}| d(n)}.
$$

Now we want to show that, for large matrices that are close to the averaging matrix  
and a suitable choice of the set $S$, the distance to the averaging can be improved with isospectral reduction; more precisely we have     
\begin{proposition}
    Given $c>0$, let $\epsilon(n) = e^{-cn} $ and $d_c(n) = \frac{1}{n} + \epsilon(n) = \frac{1}{n} + e^{-cn}$. Then for large $n$ and $|\overline{S}| = \frac{3}{4} n$ we have
    $$a_{ij} \leq d_c(n) \Rightarrow    (R_S(A,1))_{ij} \leq d_{2c}(|S|).$$
\end{proposition}
\begin{proof}
    Indeed, from the estimate above we get
    $$(R_S(A,1))_{ij} \leq \frac{d_c(n)}{1-|\overline{S}| d_c(n)} = \frac{d_c(n)}{1-\frac{3}{4}n d_c(n) }.$$
    
    We just need to show that
    $$\frac{1/n+ e^{-cn}}{1-\frac{3}{4}n (1/n + e^{-cn}) } =\frac{d_c(n)}{1-\frac{3}{4}nd_c(n)} \leq d_{2c}(|S|)= \frac{1}{n/4} + e^{-2cn/4}.$$
    The left hand side can be rewritten as
    \begin{align*}
        \frac{1/n+ e^{-cn}}{1-\frac{3}{4}n (1/n + e^{-cn})}=\frac{1/n+e^{-cn}}{1/4-3/4ne^{-cn}}=\frac{4/n+4e^{-cn}}{1-3ne^{-cn}}.
    \end{align*}
   Then what we want to show becomes
   \begin{align*}
        \frac{4/n+4e^{-cn}}{1-3ne^{-cn}}\le\frac{4}{n}+e^{-cn/2}.
   \end{align*}
   For large $n$, this is equivalent to showing
   \begin{align*}
       \frac{4}{n}+4e^{-cn}\le(1-3ne^{-cn})(\frac{4}{n}+e^{-cn/2})=\frac{4}{n}+e^{-cn/2}-12e^{-cn}-3ne^{-3cn/2},
   \end{align*}
   or
   \begin{align*}
       16e^{-cn}+3ne^{-3cn/2}\le e^{-cn/2}.
   \end{align*}
   Multiplying both sides by $e^{cn/2}$ we get
   \begin{align*}
       16e^{-cn/2}+3ne^{-cn}\le1.
   \end{align*}
   This is true for large $n$, as claimed.
\end{proof}

\subsection{Specific examples}
Here we give some examples where the matrix becomes positive after reduction. We compute the reduction and show the semi-norm and spectra before and after Isospectral reduction.
\begin{example}
Let us fix some $a \in (0, 1/2)$ and consider the set $\overline{S}$ of size two.
Let 
$$
  A_{\overline{S}\overline{S}} =
   \begin{bmatrix}
         a  &  a  \\
         a  &  a
   \end{bmatrix}.
$$
Let $|S|=n$ and
$$
 A =
\begin{bmatrix}
a & a & 1 & 1 & \cdots & 1 \\
a & a & 0 & 0 & \cdots & 0 \\
\dfrac{1-2a}{n} & 1-2a & 0 & 0 & \cdots & 0 \\
\dfrac{1-2a}{n} &  0   & 0 & 0 & \cdots & 0 \\
\vdots & \vdots & \vdots & \vdots & \ddots & \vdots \\
\dfrac{1-2a}{n} & 0 & 0 & 0 & \cdots & 0
\end{bmatrix}.
$$

In this case we have $A_{SS}$ is the zero matrix and
$$
   \tau(A) = \max{\{(1-2a)(1-\frac{1}{n}), 1-a\}} = 1-a.
$$

The reduced operator over $S$ is
$$
   R_S(A,1) = A_{S S} + A_{S\overline{S}} (I-A_{\overline{S}\overline{S}})^{-1} A_{\overline{S} S} =  A_{S\overline{S}} (I-A_{\overline{S}\overline{S}})^{-1} A_{\overline{S} S}
$$
$$
=\begin{bmatrix}
\dfrac{1-a}{n} + a & \dfrac{1-a}{n} + a & \cdots & \dfrac{1-a}{n} + a \\
\vspace{1pt}\\
\dfrac{1-a}{n} & \dfrac{1-a}{n} & \cdots & \dfrac{1-a}{n} \\
\vdots & \vdots & \ddots & \vdots \\
\dfrac{1-a}{n} & \dfrac{1-a}{n} & \cdots & \dfrac{1-a}{n}
\end{bmatrix}.
$$
The columns are all identical so $\tau(R_S(A,1)) = 0$.  The spectrum of $R_S(A,1)$ is
the set $\sigma(R_S(A,1)) = \{0,\dots,0, 1\}$, here $0$ has multiplicity $n-1$.

Let's look at the spectrum of $A$ for two specific cases.

{\bf case $n=1$:}
$$
 A =
   \begin{bmatrix}
        a & a & 1 \\
        a & a & 0 \\
        1-2a & 1-2a & 0
   \end{bmatrix}.
$$
In this case $\sigma(A) = \{2a-1, 0,  1 \}$.

{\bf case $n=2$:}
Now we get
$$
 A =
   \begin{bmatrix}
        a & a & 1 & 1 \\
        a & a & 0 & 0 \\
        \dfrac{1-2a}{2} & 1-2a & 0 & 0 \\
        \dfrac{1-2a}{2} & 0 & 0 & 0
   \end{bmatrix}.
$$
and $\sigma(A) = \{ 2a-1, 0, 0, 1 \}$, here $0$ has multiplicity $2$.
\end{example}

We can make this slightly more general.
\begin{example}
Fix $a \in (0, 1/2)$, $p \in (0, 1)$, $q=1-p$ and a column stochastic $m \times m$  matrix $B$.
Consider
$$
 A =
\begin{bmatrix}
a  &  a  & p   &  \cdots  & p   \\
a  &  a  &  0  &  \cdots  & 0 \\
\dfrac{1-2a}{m} & 1-2a &   &   &   \\
\dfrac{1-2a}{m} &  0   &   & qB   &    \\
\vdots     &   \vdots    &   &   &  \\
\frac{1-2a}{m} &  0   &   &   &
\end{bmatrix}.
$$
It's easy to show that $\tau(A) \geq  \max{ \{  q\tau(B),  \left(1-\frac{1}{m} \right)(1-2a)    \}  }$.

If we reduce $A$ over the vertices corresponding to $qB$ with value 1 for $\lambda$, we get
$$
 R_S(A,1) = q B + p
 \begin{bmatrix}
\dfrac{1-a}{m} + a  &  \dfrac{1-a}{m} + a  &  \cdots  &   \dfrac{1-a}{m} + a   \\
\dfrac{1-a}{m}       &   \dfrac{1-a}{m}     &  \cdots   &    \dfrac{1-a}{m}    \\
\vdots            &   \vdots              &   \ddots        & \vdots \\
\dfrac{1-a}{m} & \dfrac{1-a}{m}       &     \cdots  &  \dfrac{1-a}{m}
\end{bmatrix}.
$$
Since the second matrix has identical columns we have $$\tau(R_S(A,1))=q\tau(B)\le\tau(A).$$
This reduction is a positive stochastic matrix, and it is a perturbation of $qB$ when $p$ is small.
\end{example}

\begin{example}
 Still let $a \in (0, 1/2)$, $p, q >0, p+q=1$, and $B$ is an $m\times m$ column stochastic matrix. Denote by $L_i$  the average of the $i-th$ row of $B$; then we have $\sum_{i=1}^m L_i =1$. Assume that $ L_k -a/m > 0$ for any $k \in \{1, 2, \ldots, m\}$. Consider the matrix
$$
 A =
\begin{bmatrix}
a  &  a  & p &  \cdots  & p   \\
a  &  a  & 0 &  \cdots  & 0   \\
\frac{1-2a}{1-a}(L_1-a/m)   & \frac{1-2a}{m}  &   &           &    \\
\frac{1-2a}{1-a}(L_2-a/m)   & \frac{1-2a}{m}  &   & qB      &    \\
\vdots & \vdots  &   &           &     \\
\frac{1-2a}{1-a}(L_m - a/m) &  \frac{1-2a}{m}  &   &           &
\end{bmatrix}.
$$
We can show that
$$
 \tau(A) \geq \max{ \{  q \tau(B), \frac{1-2a}{2(1-a)}( |L_1-\frac{1}{m}| + \cdots + |L_N - \frac{1}{m} | ) \}  }.
$$

Again we reduce $A$ over the vertices corresponding to $qB$ with value 1 for $\lambda$,
$$
 R_S(A,1) = q B + p
 \begin{bmatrix}
L_1  &    L_1  &  \cdots  &   L_1   \\
L_2  &    L_2  &  \cdots  &   L_2    \\
\vdots &  \vdots & \ddots & \vdots \\
L_m &     L_m  &  \cdots  &   L_m
 \end{bmatrix}.
$$
We have $\tau(R_S(A,1))= q \tau(B)\le\tau(A)$ and the reduction is a positive stochastic matrix.
\end{example}

Notice that the matrix
$$
   L =
\begin{bmatrix}
L_1  &    L_1  &  \cdots  &   L_1   \\
L_2  &    L_2  &  \cdots  &   L_2    \\
\vdots &  \vdots & \ddots  & \vdots \\
L_m &     L_m  &  \cdots  &   L_m
\end{bmatrix}
$$
is column stochastic. It has eigenvalue $0$ with corresponding eigenspace $[1, 1, \ldots, 1]^{\perp}$  and eigenvalue $1$, whose eigenvector is $[L_1, L_2, \ldots, L_m]$.

 Although $L_i$'s correspond to the averages of the rows of $B$,
 they can be seen as free parameters; in particular, we can consider the choice $L_i = 1/m$
 for all $i \in \{1, 2, \ldots, m\}$. In this case the reduced matrix is
 $$
 R_S(A,1) = q B + (1-q)
\begin{bmatrix}
1/m  &    1/m  &  \cdots  &   1/m  \\
1/m  &    1/m  &  \cdots  &   1/m  \\
\vdots & \vdots & \ddots & \vdots \\
1/m &     1/m  &  \cdots  &   1/m
\end{bmatrix},
$$
the famous Google matrix. Here $q$ is the damping factor.

\begin{example}
Adopting all the same notations as the previous example, we now consider the case where the set $\overline{S}$ has only one point.
$$
 A =
\begin{bmatrix}
a  &  p  &  \cdots  & p   \\
(L_1-a/m)   &   &           &    \\
(L_2-a/m)   &  & qB  &    \\
\vdots          &   &           &     \\
(L_m - a/m)          &   &           &
\end{bmatrix},
$$
then
$$
\tau(A)=\max{ \left\{  q \tau(B),  \max_j{ \{ \frac{1}{2}[ |a-p| + |L_1-a/m-qb_{1j} | + \ldots +|L_m -a/m -qb_{mj} ]\} }\right\} }.
$$
When $q$ is close to zero $\tau(A)$ approaches $1-a$.

For the reduced operator we have
$$
  R_S(A,1) = q B +  \frac{p}{1-a}  
\begin{bmatrix}
L_1-a/m  &    L_1-a/m  &  \cdots  &   L_1-a/m  \\
L_2-a/m  &    L_2-a/m  &  \cdots  &   L_2-a/m  \\
\vdots & \vdots  & \ddots & \vdots \\
L_m-a/m &     L_m-a/m  &  \cdots  &   L_m-a/m 
\end{bmatrix}.
$$
Again the reduction is a positive stochastic matrix and $\tau(R_S(A,1))=q \tau(B)\le\tau(A)$.
\end{example}

\begin{example}
 For the last example let's consider an $n \times n$ column stochastic
 matrix $A$  such that, for some $m < n$,
 the entry $a_{ik}$  is positive if and only if  $k \in \{i - (m-1), \ldots, i + (m-1) \}$. Take $S = \{1, 2, \ldots, m\}$. In this case, the submatrix $A_{SS}$ is positive, and so
 $$
   R_S(A,1) = A_{SS} + A_{S\overline{S}}(I-A_{\overline{S} \overline{S}})^{-1} A_{\overline{S} S}
 $$
 is also positive.

 Observe that $A$ is irreducible so that $(I-A_{\overline{SS}})^{-1}=I+A_{\overline{SS}}+A^2_{\overline{SS}}+\dots$ is positive.
\end{example}

\section{Discussion and Future Work}

In Section \ref{spec}, we proved that the diameter semi-norm always decreases after isospectral reductions. Our numerical experiments in Section \ref{comp} indicate that the inner spectral radius tends to decrease with high probability (see Figure \ref{8-10}). Additionally, in Section \ref{sec-positivematrices}, we provided examples of nonnegative matrices that turn positive after isospectral reductions. Notice that the diameter semi-norm tends to decrease as the entries of the matrix become more positive and uniform. While the semi-norm serves as an upper bound for the inner spectral radius and always decreases after reduction, the inner spectral radius does not consistently follow suit. An example illustrating this is provided below.

\begin{example}For the stochastic matrix
    $$A=\begin{bmatrix}
        0 & 0.9 & 0 \\
        0.5 & 0 & 1 \\
        0.5 & 0.1 & 0
    \end{bmatrix},\tau(A)=1,\rho_i(A)=0.6708.$$ Let's reduce $A$ to the first two nodes and plug in value $1$ for $\lambda$, we get
    $$R_S(A,1)=\begin{bmatrix}
        0 & 0.9 \\
        1 & 0.1
    \end{bmatrix},\tau(R_S(A,1))=0.9,\rho_i(R_S(A,1))=0.9.$$
    Here after reduction $R$ is still a stochastic matrix, its semi-norm is lower than $A$ and its inner spectrum is higher than $A$.
\label{counter--}
\end{example}

We are interested in the following questions for future work.

\begin{itemize}
    \item Can we prove that the inner spectral radius decreases with high probability
 after isospectral reductions for specific classes of stochastic matrices? For example, stochastic matrices with unbounded first moment or stochastic matrices associated with digraphs generated by the Erdos-Rényi or Barabási-Albert models.
 \item For these specific classes of stochastic matrices, can we estimate the factor of decrease?
 \item Are there simple checkable conditions for the inner spectral radius to decrease after isospectral reductions? For example, can we give simple heuristics to choose subsets $S\subset\{1,...,n\}$ of a given size  which ensures the decrease of the inner spectral radius?
 \item Can we describe simple heuristics to remove nodes one by one until reaching a subset $S$ such that $A_{\overline{SS}}$ is well conditioned and the inner spectral radius of $R_S(A,1)$ gets below a given threshold?
\end{itemize}
Such conditions could lead to faster algorithms based on isospectral reductions. 
\section*{Acknowledgments}
PD was supported  by Funda\c{c}\~{a}o para a Ci\^{e}ncia e a Tecnologia, through the project UID/MAT/04561/2013.
MJT was partially financed by Portuguese Funds through FCT  - `Funda\c{c}\~ao para a Ci\^encia e a Tecnologia' within the Projects UIDB/00013/2020 (DOI: 10.54499/UIDB/00013/2020)  and UIDP/00013/2020  (DOI: 10.54499/UIDP/00013/2020).
AB, PD, LS and MJT were partially supported by the Project
"New trends in Lyapunov exponents"(PTDC/MAT-PUR/29126/2017).

\bigskip



\bibliographystyle{elsarticle-num} 


\end{document}